\documentclass{amsart}
\usepackage{amsmath}
\usepackage{amssymb}
\usepackage[curve,arrow,matrix]{xy}
\theoremstyle{plain}
\newtheorem{thm}{Theorem}[section]
\newtheorem{prop}[thm]{Proposition}
\newtheorem{lemma}[thm]{Lemma}
\newtheorem{cor}[thm]{Corollary}

\newcommand{\deqno}{\refstepcounter{thm}(\thethm)}

\theoremstyle{definition}
\newtheorem{defn}[thm]{Definition}
\newtheorem{remk}[thm]{Remark}

\newcounter{item}
\newenvironment{rlist}[1][1]{\begin{list}
  {\textup{(\roman{item})}}{\usecounter{item} \setcounter{item}{#1}\addtocounter{item}{-1}
  \setlength{\itemsep}{0ex}
  \setlength{\topsep}{0ex} \setlength{\parsep}{0ex} \setlength{\labelwidth}{15mm}
   \setlength{\leftmargin}{10mm} } }{\end{list}}

\newcommand{\ba}{{\bf a}}
\newcommand{\bb}{{\bf b}}

\newcommand{\bx}{{\bf x}}
\newcommand{\by}{{\bf y}}
\newcommand{\bs}{{\bf s}}

\newcommand{\Fe}{{\bf F}^e}
\newcommand{\nn}{\mathbb{N}}
\newcommand{\qq}{\mathbb{Q}}
\newcommand{\zz}{\mathbb{Z}}

\newcommand{\tensor}{\otimes}
\newcommand{\Hom}{\textup{Hom}}

\newcommand{\Proj}{\textup{Proj}}

\newcommand{\im}{\textup{Im}}

\newcommand{\gr}{\textup{gr}}
\newcommand{\GR}{\textup{GR}}

\newcommand{\supp}{\textup{supp}}

\newcommand{\ra}{\rightarrow}

\newcommand{\incl}{\hookrightarrow}

\newcommand{\isom}{\cong}

\begin{document}

%\begin{frontmatter}
\title{Tight Closure of Finite Length Modules in Graded Rings}
\author{Geoffrey D. Dietz}
\thanks{The author was supported in part by a VIGRE grant from
 the National Science Foundation.}
\thanks{Journal of Algebra, 306 (2006), No. 2, 520--534.}
%\corauth[cor]{Corresponding Author}
\email{gdietz@member.ams.org}
\address{Department of Mathematics, Gannon University, Erie, Pennsylvania
16541}

\subjclass[2000]{Primary 13A35; Secondary 13B22, 13C11}

\keywords{Tight closure, plus closure, injective hulls}

%    General info
\begin{abstract}
In this article, we look at
how the equivalence of tight closure and plus closure (or Frobenius
closure) in the homogeneous $m$-coprimary case
implies the same closure equivalence in the non-homogeneous $m$-coprimary
case in standard graded rings.
Although our result does not depend upon dimension,
the primary application is based on results known in dimension 2
due to the recent work of H.\ Brenner. We also demonstrate
a connection between tight closure and the $m$-adic closure
of modules extended to $R^+$ or $R^\infty$.
We finally show that unlike the Noetherian case, the injective
hull of the residue field over $R^+$ or $R^{\infty}$
contains elements that are not killed by any power of the maximal
ideal of $R$. This fact presents an obstruction to
one possible method of extending our main result to all modules.
\end{abstract}

\maketitle

\section{Introduction}
\label{intro}

Tight closure theory, since its introduction
by M.\ Hochster and C.\ Huneke in the late 1980s, has been
an important method of working in positive characteristic.
The \textit{tight closure} of an ideal is defined to be
$$
I^* := \{x\in R\,|\, \exists c\in R^\circ \textup{ such that }
cx^{p^e}\in I^{[p^e]}\quad \forall e\gg 0\}
$$
where $p$ is the characteristic of the ring, $R^\circ$ is
the complement in $R$ of the minimal primes
of $R$, and $I^{[p^e]}$ is the ideal generated by all $p^e$th
powers of the elements of $I$.
Some of the success of tight closure has been due to its ability to
tie together ideas that were previously not known to be connected,
generalize theorems, and simplify proofs.
Some of the many examples include generalizations of
the Brian\c{c}on-Skoda Theorem and the Hochster-Roberts Theorem.
See \cite{HH90} for a general introduction.

Deciding whether tight closure computations
commute with localization or not
remains a very elusive goal. See \cite{AHH} for results on this problem.
Closely related to the localization
issue is the question of whether the tight closure of an ideal $I$
in a positive characteristic domain $R$
is simply the contracted-expansion of $I$ to $R^+$, the integral closure
of $R$ in an algebraic closure of its fraction field, i.e.,
does $I^*=IR^+\cap R$? (The ideal $IR^+\cap R$ is called the
\textit{plus closure} of $I$.) Since the plus closure operation
can easily be shown to commute with localization, one could solve the
localization problem by proving the above equality.

In \cite{Sm94}, K.E.\ Smith made a tremendous contribution
to this problem by proving that $I^*=IR^+\cap R$ for ideals
generated by partial systems of parameters in excellent
local domains of positive characteristic. Smith
also showed that tight closure and plus closure are equal
for ideals generated by part of a homogeneous system
of parameters in an $\nn$-graded domain of
positive characteristic. See \cite{Sm95} or Theorem
\ref{param+tight} here.

Recently H.\ Brenner made a major breakthrough on this
problem when he showed that tight closure and plus closure are
equivalent for homogeneous ideals in certain 2-dimensional
graded rings. See \cite{Br05} and \cite{Br3} and Theorems
\ref{Br-ellcurve}, \ref{Br-curve}, and \ref{Br-ellFr} here.
Brenner's work relies on
a correspondence between tight closure and
cohomological behavior of bundles on projective
varieties (see \cite{Br03}).
Although it is not made explicit in his work, the same methods
show the equivalence for finitely generated graded modules
over the same class of rings. See Section \ref{TCconseq}
for details.

Inspired by Brenner's work, we have studied how
one can obtain an equivalence of tight closure and plus closure
for more general ideals and modules given that one
has the equivalence for homogeneous ideals and modules.
Our main result is the following theorem. \vspace{3mm}

\noindent
\textbf{Theorem \ref{fltc=+cl}.} \textit{Let $(R,m)$ be a standard graded
$K$-algebra (see Section \ref{mainsection})
of characteristic $p>0$. Suppose that $R$ is
a domain, and $K$ is algebraically closed.
If $N_M^* = N_M^+ = N_M^{+\gr}$
for all finitely generated graded $R$-modules
$N\subseteq
M$ such that $M/N$ is $m$-coprimary, then the same is true for \textup{all}
finitely generated modules $N\subseteq M$ such that $M/N$ is $m$-coprimary.}
\vspace{3mm}

As a result,
we can apply our theorem to the cases where Brenner's work
is valid to increase the class of ideals and modules where
tight closure equals plus closure.

Unlike the work of Brenner,
our methods are entirely algebraic and rely on
injective modules over a graded subring, $R^{+\GR}$,
of $R^+$. This led us to the study of injective hulls over
$R^{+\GR}$ and $R^\infty$ in an attempt to extend our
result beyond the $m$-coprimary case. We present
a submodule of the injective hull
$E_{A^\infty}(K)$, where $A$ is either a polynomial
ring or a formal power series ring, that
we use to show that the
injective hulls $E_{R^{+\GR}}(K)$, $E_{R^+}(K)$, and $E_{R^\infty}(K)$
behave far differently from the Noetherian case, as
these modules contain elements that are not killed
by any power of the maximal ideal of $R$.
See Theorems \ref{badperfhull} and \ref{badinthull}. As a consequence
we have not been able to use these injective modules to extend
Theorem \ref{fltc=+cl} to general modules.

\section{Notation and Background}

Before we state our results, we
survey the theory that forms the
foundation and provides the motivation for our work.
All rings throughout are commutative with identity and are Noetherian
unless noted otherwise. All modules are unital.

\subsection{The Frobenius Endomorphism}

We will always let $p$
denote a positive prime number, and $q$ will denote $p^e$, a power
of $p$.
Every characteristic $p$ ring $R$ comes equipped with a \textit{Frobenius
endomorphism} $F_R:R\ra R$, which maps
$r\mapsto r^p$. Composing this map with itself we obtain
$F_R^e:R\ra R$,
which map $r\mapsto r^q$. Closely associated to these maps
are the \textit{Peskine-Szpiro} (or \textit{Frobenius})
\textit{functors} $\Fe_R$.
If we let $S$ denote the ring
$R$ viewed as an $R$-module via the $e^{th}$-iterated Frobenius
endomorphism, then $\Fe_R$
is the covariant functor $S\tensor_R -$ which
takes $R$-modules to $S$-modules and so takes
$R$-modules to $R$-modules since $S=R$ as a ring. Specifically,
if $R^m\ra R^n$ is a map of free $R$-modules given by the matrix
$(r_{ij})$, then we may apply $\Fe_R$ to this map to obtain a map
between the same $R$-modules given by the matrix
$(r^q_{ij})$. For cyclic modules
$R/I$, $\Fe_R(R/I)=R/I^{[q]}$, where
$$
I^{[q]} := (a^q \,|\, a\in I)R
$$
is the \textit{$q^{th}$ Frobenius power} of the ideal $I$. If the ideal $I$
is finitely generated, then $I^{[q]}$ is also the ideal generated
by the $q^{th}$ powers of a finite generating set for $I$.
In a similar manner,
for modules $N\subseteq M$,
$$
N_M^{[q]} := \im(\Fe_R(N)\ra \Fe_R(M)),
$$
and we will denote
the image of $u\in N$ inside of $N_M^{[q]}$ by $u^q$.

When $R$ is reduced, we define
$R^{1/q}$ to be
the ring obtained by adjoining to $R$ all $q^{th}$ roots of elements
in $R$. In this setting, the inclusion $R\incl R^{1/q}$ is isomorphic
to the inclusion $F^e_R:R\incl R$, identifying $R^{1/q}$ with $R$
via the isomorphism $r^{1/q}\mapsto r$.
Therefore the Peskine-Szpiro
functor $\Fe_R$ is isomorphic to $R^{1/q}\tensor_R -$
after identifying $R^{1/q}$ with $R$,
and the Frobenius
power $I^{[q]}$ can be identified in the same sense
with the extension $IR^{1/q}$.
We denote by $R^\infty$ the \textit{perfect closure} of $R$.
The ring $R^\infty$ is constructed by adjoining all $q^{th}$
roots to $R$, for all $q$.
In general, $R^\infty$ is not a Noetherian ring.

\begin{defn} For a Noetherian ring $R$ of positive characteristic $p$ and
finitely generated $R$-modules
$N\subseteq M$, the \textit{Frobenius closure} of $N$ in $M$ is the
submodule
$$
N_M^F := \{ u\in M \,|\, u^q\in N^{[q]}_M,\, \mbox{for some}\, q\}.
$$
\end{defn}

\begin{lemma}\label{fclosure} If $R$ is a reduced Noetherian ring
of characteristic $p>0$, then for
finitely generated $R$-modules $N\subseteq M$, the following are
equivalent:
\begin{rlist}
\item $u\in N_M^F$.
\item $1\tensor u \in \im(R^{1/q}\tensor N \ra R^{1/q}\tensor M)$,
for some $q$.
\item $1\tensor u \in \im(R^\infty\tensor N \ra R^\infty\tensor M)$.
\end{rlist}
In the case of ideals, $I^F = \bigcup_q IR^{1/q}\cap R = IR^\infty \cap R$.
\end{lemma}

\subsection{Tight Closure}

Let the complement in $R$ of the set of minimal primes be
denoted by $R^\circ$.

\begin{defn} For a Noetherian ring $R$ of characteristic $p>0$
and finitely generated modules $N\subseteq M$, the \textit{tight
closure} $N_M^*$ of $N$ in $M$ is
$$
N_M^* := \{ u\in M\,|\, cu^q\in N_M^{[q]}\, \mbox{for all}\, q\gg 1, \,
\mbox{for some}\, c\in R^\circ\}.
$$
In the case that $M=R$ and $N=I$, $u\in I^*$
if and only if
there exists $c\in R^\circ$ such that $cu^q\in I^{[q]}$, for all
$q\gg 1$.
\end{defn}

Using the existence of \textit{test elements} (see \cite{HH90,HH94sm}), 
we have some useful
characterizations of tight closure using $R^{1/q}$ and
$R^\infty$.

\begin{lemma}\label{tcchar} Let $N\subseteq M$ be finitely
generated $R$-modules, where $R$ is reduced of positive
characteristic $p$ and has a test element. Then the following are equivalent:
\begin{rlist}
\item $u\in N_M^*$.
\item $c^{1/q}\tensor u \in \im(R^{1/q}\tensor N \ra R^{1/q}\tensor
M)$, for all $q\gg 1$ and some $c\in R^\circ$.
\item $c^{1/q}\tensor u \in \im(R^\infty\tensor N \ra R^\infty\tensor
M)$, for some (or every) test element $c$ and all $q\geq 1$.
\end{rlist}
\end{lemma}

\subsection{Plus Closure}
\label{intro+cl}

Let $R$ be a domain, and let $R^+$ denote the integral
closure of $R$ in an algebraic closure of its
fraction field. $R^+$ is called the
\textit{absolute integral closure of $R$} and is not
Noetherian in general. This is an important
ring due to the remarkable result of Hochster and
Huneke.

\begin{thm}[Theorem 5.15, \cite{HH92}]\label{R+bigCM}
Let $R$ be an excellent local
domain of positive characteristic. Then $R^+$ is a big Cohen-Macaulay
$R$-algebra, i.e., every system of parameters of $R$ is a regular sequence
on $R^+$.
\end{thm}

Hochster and Huneke also provide a graded version of the above theorem.
If $R$ is an $\nn$-graded domain, then $R^{+\GR}$
denotes a maximal
direct limit of module-finite, $\qq_{\geq 0}$-graded extension
domains of $R$. For the construction and properties of this ring, see
\cite[Lemma 4.1]{HH92}. There is also an
$\nn$-graded direct summand of $R^{+\GR}$, which is
denoted $R^{+\gr}$. Neither of these
rings is Noetherian in general.

\begin{thm}[Theorem 5.15, \cite{HH92}]\label{grR+bigCM}
If $R$ is an $\nn$-graded domain of positive characteristic
with $R_0=K$ and $R$ a finitely
generated $K$-algebra, then $R^{+\GR}$ and
$R^{+\gr}$ are both graded big Cohen-Macaulay $R$-algebras
in the sense that every homogeneous system of parameters
of $R$ is a regular sequence on $R^{+\GR}$ and $R^{+\gr}$.
\end{thm}

There is a closure operation associated to each
of the rings $R^+$, $R^{+\GR}$, and $R^{+\gr}$.
Since $R^{+\gr}$ is a direct summand of $R^{+\GR}$ as an $R^{+\gr}$-module,
the latter two rings yield equivalent closure operations.

\begin{defn}\label{+cldefn} Given an excellent, local
(resp., Noetherian, $\nn$-graded) domain $R$ of positive characteristic, 
let $S=R^+$ (resp., $S=R^{+\GR}$ or $S=R^{+\gr}$), and let
$N\subseteq M$ be $R$-modules. The \textit{plus closure} $N_M^+$
(resp., \textit{graded-plus closure} $N_M^{+\gr}$)
of $N$ in $M$ is
$$
\{u\in M \,|\, 1\tensor u\in \im(S\tensor_R N\ra
S\tensor_R M)\}.
$$
If $M=R$ and $N=I$, then $I^+$ (resp., $I^{+\gr}$)
equals $IS\cap R$.
\end{defn}

It is straightforward to show that $N_M^F\subseteq N_M^+\subseteq N_M^*$,
for all finitely generated modules, and that
$N_M^F\subseteq N_M^{+\gr}\subseteq N_M^+\subseteq N_M^*$ in the
graded case. As mentioned earlier,
it is possible that tight closure in positive characteristic is just plus closure.
Some of the most important
results in this direction come from K.E.\ Smith.

\begin{thm}[\cite{Sm94}, \cite{Sm95}]\label{param+tight} Let $R$ be an excellent local
(resp., $\nn$-graded) Noetherian domain of characteristic $p>0$ (with $R_0$ a field), 
and let $I=(x_1,\ldots,x_k)$ be an ideal generated by
part of a (homogeneous) system of parameters. Then $I^*= I^+$ (resp.,
$I^* = I^{+\gr}= I^+$).
\end{thm}

Recently, H.\ Brenner has developed results in dimension 2
that show that tight closure and graded-plus closure
are equivalent for homogeneous ideals in certain graded rings.

\begin{thm}[Theorem 4.3, \cite{Br05}]\label{Br-ellcurve}
Let $K$ be an algebraically closed field of positive characteristic,
and let $R$ be the homogeneous coordinate ring of an elliptic curve
(i.e., $R$ is a standard graded normal $K$-algebra of dimension
2 with $\dim_K [H_m^2(R)]_0 = 1$, where $m$ is the
homogeneous maximal ideal of $R$).
Let $I$ be an $m$-primary graded
ideal in $R$. Then $I^{+\gr}=I^+=I^*$.
\end{thm}

For example, the result above applies when
$R=K[x,y,z]/(F)$ is normal, where $F$ is homogeneous of degree 3.

\begin{thm}[Theorem 4.2, \cite{Br3}]\label{Br-curve} Let $K$ be the
algebraic closure of a finite field. Let $R$ denote an $\nn$-graded
2-dimensional domain of finite type over $K$. Then for every homogeneous
ideal $I$, we have $I^{+\gr}=I^+=I^*$.
\end{thm}

In the case of an elliptic curve with Hasse invariant 0
(see \cite[pp.\ 332--335]{Ha}),
Brenner showed that tight closure is the same as Frobenius closure.

\begin{thm}[Remark 4.4, \cite{Br05}]\label{Br-ellFr}
If $R$ is the homogeneous coordinate ring of an
elliptic curve of positive characteristic $p$
with Hasse invariant 0 defined over an
algebraically closed field, then $I^* = I^F$ for
all $m$-primary graded ideals of $R$, where
$m$ is the homogeneous maximal ideal
of $R$.
\end{thm}

\section{New Cases Where Tight Closure is Plus Closure}
\label{mainsection}

Before proving our main result, Theorem \ref{fltc=+cl},
we need to establish
some lemmas and notation. If $S$ is any $\qq$-graded (not necessarily
Noetherian) ring, then for any $n\in\qq$ let $S_{\geq n} =
\bigoplus_{i\geq n} S_i$.
Similarly define $S_{>n}$. We will say that
an $\nn$-graded ring $R$ is a \textit{standard graded $R_0$-algebra} if
$R$ is finitely generated over $R_0$ by elements of degree 1.
For the rest of the section, let
$m=\bigoplus_{i>0} R_i$, the homogeneous maximal ideal of $R$.

\begin{lemma}\label{infbnd} Let $R$ be a reduced standard graded $K$-algebra
of positive characteristic $p$, and let $S=R^\infty$. Then
there exists $c\in\nn$ such that $S_{\geq n+c}\subseteq
m^nS$ for any $n\geq 1$. As a consequence, $[S/m^nS]_j = 0$ for all
$j\gg 0$.
\end{lemma}
\begin{proof}
Let $m$ be generated by $x_1,\ldots,x_\mu$, each of degree 1.
Put $c=\mu-1$ (if $\mu = 0$,
i.e., $R=K$, put $c=0$). Since
$S_{>0}=\cup_q m^{1/q}$, if $f\in S$ is homogeneous of degree
at least $n+c$, then
$f$ is a sum of terms $dx_1^{\alpha_1}\cdots x_\mu^{\alpha_\mu}$
such that $d\in S_0$ and $\sum \alpha_i\geq n+c$. If we write
$\alpha_i = [\alpha_i] + r_i$, where $0\leq r_i<1$ for all $i$, then
$$\sum [\alpha_i] = \sum \alpha_i -\sum r_i \geq n+c - \sum r_i >
n+c-\mu.$$
Therefore, $\sum [\alpha_i] \geq n+c-\mu +1 = n$, and so
$f\in m^nS$. The second claim now follows since $j\geq n+c$. \end{proof}

To prove a similar result for $R^{+\GR}$,
we will need a graded-plus closure version of the
Brian\c{c}on-Skoda Theorem. The
original tight closure generalization can be found in
\cite[Theorem 5.4]{HH90}. Hochster and Huneke also strengthened this
result to a version for plus closure in \cite[Theorem 7.1]{HH95}.
We will adapt their proof to obtain a graded-plus closure version of
the Brian\c{c}on-Skoda Theorem.

\begin{thm}\label{BS+GR} Let $R$ be a positively graded Noetherian
domain of positive characteristic. Let $I$ be a homogeneous
ideal generated by at most $d$ homogeneous elements, let $k\in\nn$,
and let $u\in \overline{I^{d+k}}$ with $u$ homogeneous. Then
$u\in I^{k+1}S\cap R$, where $S=R^{+\GR}$ or $S=R^{+\gr}$.
\end{thm}
\begin{proof}
We will first use $R$, $I$, and $u$ to construct a triple $(A,J,v)$ and a
degree-preserving map to the triple $(R,I,u)$ such that
$A$ is a positively graded Noetherian
domain, $J$ is a homogeneous ideal of $A$ with
at most $d$ generators, $v\in \overline{J^{d+k}}$
is homogeneous in $A$, 
$I=JR$, and $v\mapsto u$. We will directly prove  
the theorem holds for $(A,J,v)$ and then show that this case implies
the result for $(R,I,u)$.

Since $u\in R$ is integral over $I^{d+k}$ and homogeneous,
$u$ satisfies a homogeneous monic polynomial
$$
z^n+r_1z^{n-1}+\cdots + r_n=0,
$$
where $\deg z=\deg u$, $\deg r_j = j\deg u$, $r_j\in (I^{d+k})^j$, and
(without loss of generality) $r_n\neq 0$.
Each $r_j$ can be written
as a homogeneous $R$-linear combination of
monomials $a_1^{\nu_1}\cdots a_d^{\nu_d}$
in the generators $a_1,\ldots,a_d$ of $I$,
where $\nu_1+\cdots +\nu_d = (d+k)j$. Thus, the coefficient
of the monomial $a_1^{\nu_1}\cdots a_d^{\nu_d}$ is zero or has
$$
\mbox{degree}\, = \deg r_j -(\nu_1\deg a_1 +\cdots + \nu_d\deg a_d)
$$
since $R$
is positively graded. Without
loss of generality, we may order the generators of $I$ so that
$\deg a_1\leq \cdots\leq\deg a_d$. Then $r_n\neq 0$ implies that
$\deg a_1^{(d+k)n}\leq \deg r_n$. If not, then
$\deg r_n<\nu_1\deg a_1+\cdots +\nu_d \deg a_d$, for all $\nu_i$
such that $\nu_1+\cdots +\nu_d=(d+k)n$, and so the coefficient of every
monomial in the expansion of $r_n$ must be zero, a contradiction.

Let $x_1,\ldots,x_d$ be indeterminates over $K=\zz/p\zz$ with $\deg x_i=\deg a_i$.
For every monomial $\mu=x_1^{\nu_1}\cdots x_d^{\nu_d}$, where
$\nu_1+\cdots +\nu_d = (d+k)j$ for $1\leq j\leq n$, let $y_\mu$ be
an indeterminate with $\deg y_\mu = \deg r_j - (\nu_1\deg x_1 +
\cdots +\nu_d\deg x_d)$. Let
$$F(\bx,\by,z) = z^n + \sum_{j=1}^n
\left(\sum_{\mu\in C_j} y_\mu \mu\right)z^{n-j},$$
where $\bx =x_1,\ldots,x_d$, $\by = \{y_\mu\, |\, \deg y_\mu\geq 0\}$,
and
$$C_j=\{\mu=x_1^{\nu_1}\cdots x_d^{\nu_d}\, |\, \nu_1+\cdots +\nu_d = (d+k)j\}.$$
Then $F$ is homogeneous of degree $n\deg z=n\deg u$ as
$$\deg (y_\mu \mu)z^{n-j} = \deg r_j +(n-j)\deg z = j\deg u
+n\deg z -j\deg z = n\deg z.$$

Therefore, $K[\bx,\by,z]$ is a positively graded Noetherian ring,
and the homomorphism $K[\bx,\by]\ra R$, given by $x_i\mapsto a_i$ and mapping
$y_\mu$
to the coefficient of $a_1^{\nu_1}\cdots a_d^{\nu_d}$, is degree-preserving. 
Moreover, the composite map $K[\bx,\by,z]\ra R[z]\ra R$, where
$z\mapsto u$, sends $F(\bx,\by,z)\mapsto z^{n}+r_{1}z^{n-1}+\cdots
+r_{n}\mapsto 0$. (Since $R$ is positively graded, the coefficient
of $a_{1}^{\nu_{1}}\cdots a_{d}^{\nu_{d}}$ is $0$ if
$\deg r_{j}<\nu_{1}\deg a_{1}+\cdots +\nu_{d}\deg a_{d}$ so
that we did not need a $y_{\mu}\mu$ term in $F$ when $\deg y_{\mu}<0$.)

Put $A:=K[\bx,\by,z]/F(\bx,\by,z)$, $J:=(\bx)A$, and $v:=z$ in $A$. Then $A$ is a
positively graded Noetherian ring of positive characteristic, $J$ is
a homogeneous ideal generated by at most $d$ homogeneous elements, and
$v$ is homogeneous and in $\overline{J^{d+k}}$. It is clear
from the construction of $A$ that $A\ra R$ is a degree-preserving map,
that $I=JR$, and that $v\mapsto u$ under the map.
To see that $A$ is also a domain, we will show that $F$ is irreducible.
Indeed, let $N=(d+k)n$, and let $\mu$ be the monomial $x_{1}^{N}$ that
occurs when $j=n$ in the summation for $F$.
As we noted earlier, $r_{n}\neq 0$ implies that
$$\deg y_{\mu} = \deg r_{n} - N\deg x_{1} = \deg r_{n} -N\deg a_{1}\geq
0,$$
and so $F$ is linear in $y_{\mu}$ with coefficient $x_{1}^{N}$ for
$y_{\mu}$ and a relatively prime constant term containing $z^{n}$.

We next show that the theorem holds for the triple
$(A,J,v)$ constructed above. As $A^{+\gr}$ is a direct
summand of $A^{+\GR}$, it is enough to show that
$v\in JA^{+\GR}\cap A$. Since
$A$ is a positively graded, finitely generated $K$-algebra, we may
regrade if necessary so that it is $\nn$-graded without changing
$A^{+\GR}$. Since $A/(\bx)\isom K[\by,z]/z^{n}$, the sequence
$x_{1},\ldots,x_{d}$ forms part of a homogeneous system of parameters.
We can now apply Theorem \ref{param+tight} to the ring $A$ and ideal $J$ to
see that $J^{*}= JA^{+\GR}\cap A$. By
the Generalized Brian\c{c}on-Skoda Theorem (\cite[Theorem 5.4]{HH90}), 
the theorem holds for $(A,J,v)$.

We finally claim that the theorem holds for the original triple
$(R,I,u)$. (We again only need to show the $R^{+\GR}$ case.) Since the map
$A\ra R$ extends to $A^+\ra R^+$, we can restrict
this map to obtain $A^{+\GR}\ra R^{+\GR}$ (The homogeneous monic equation
satisfied by an element $a$ of $A^{+\GR}$ maps to a homogeneous
monic equation over $R$ satisfied by the image of $a$.)
Therefore, $u\in I^{k+1}R^{+\GR}$ as $v\in J^{k+1}A^{+\GR}$, 
$v\mapsto u$, and $JR = I$.
\end{proof}

\begin{lemma}\label{gradedbnd} Let $R$ be a standard graded $K$-algebra
domain of characteristic $p>0$, and let $S=R^{+\GR}$ or $S=R^{+\gr}$.
Then there exists $c\in\nn$
such that $S_{\geq n+c}\subseteq m^nS$, for any $n\geq 1$. Moreover,
$[S/m^nS]_j = 0$ for all $j\gg0$. \end{lemma}
\begin{proof}
Let $m$ be generated
by $\mu$ elements. Let $c=\mu-1$ (if $\mu = 0$, let $c=0$), and let $f\in S$
be homogeneous
of degree $D\geq n+c$. Then $f$ satisfies a monic polynomial equation
$f^t+r_1f^{t-1}+\cdots + r_t = 0$ such that $r_i$ is homogeneous of
degree $iD$ in $R$ or $r_i = 0$ if $iD\not\in\nn$. Therefore, $r_i\in
m^{i(n+c)} = (m^{n+c})^i$ for all $i$ as $m$ is generated in degree 1.
Since $f\in S$, there exists a positively graded
module-finite extension domain $T$ of $R$ such that $f\in T$. Thus,
$f\in \overline{(mT)^{n+c}} = \overline{(mT)^{\mu+n-1}}$. By
Theorem \ref{BS+GR}, $f\in m^nT^{+\GR}$, but $T^{+\GR}=R^{+\GR}$,
and so $f\in m^n R^{+\GR}$. Since $R^{+\gr}$ is a direct summand of
$R^{+\GR}$, we also have $f\in m^n R^{+\gr}$. The second
claim follows using $j\geq n+c$.  \end{proof}

Our main result will depend upon showing that
$\Hom_K(S/m^n S, K)$ is $\zz$-graded as an $R$-module
when $S$ is $R^{+\GR}$, $R^{+\gr}$, or $R^{\infty}$.

\begin{prop}\label{grading} Let $R$ be a standard graded
$K$-algebra of characteristic $p>0$. Suppose $R$
is reduced
(respectively, a domain). Let $S=R^\infty$ (resp., $S=R^{+\GR}$ or
$S=R^{+\gr}$). Then
for any $n\geq 1$, $\Hom_K(S/m^n S, K)$ is a $\zz$-graded $R$-module.
\end{prop}
\begin{proof} $S$ has a natural $\nn[1/p]$-grading (resp.,
$\qq_{\geq 0}$-grading or $\nn$-grading) induced by
the grading on $R$. Thus, $S/m^nS$ is also graded as $m^nS$ is a
homogeneous ideal. Let $W_j$ be the $K$-span of all
homogeneous elements of degree $\delta$ such that $j-1<\delta\leq j$.
This gives $S/m^nS$ an $\nn$-grading as an $R$-module, where
$W_j = 0$ for all $j<0$ and $j\gg 0$ by Lemma \ref{infbnd}
(resp., Lemma \ref{gradedbnd}).

In $\Hom_K(S/m^n S, K)$, let $V_{-j}$ be the $K$-span of all
functionals $\phi$ such that $\phi(W_i)$ is not 0 if and only if $i =
j$.
If $r\in R$ is homogeneous of degree $d$, and $\phi\in V_{-j}$, then
$r\phi(W_i)=\phi(rW_i)\subseteq \phi(W_{i+d})$ which is nonzero if
and
only if $i = -d+j$. Thus, $R_dV_{-j}\subseteq V_{-j+d}$. It is
clear that the intersection of any $V_{-j}$ with the sum of the
others is trivial and that $\sum_j V_{-j} \subseteq \Hom_K(S/m^n S,
K)$. Now, if $\psi\in \Hom_K(S/m^n S, K)$, and $s$ has homogeneous
components $s_i$, then let $\psi_{-j}(s)=\psi(s_j)$ so that
$\psi_{-j}\in V_{-j}$. Then $\psi
=\sum_j \psi_{-j}$, where the sum is finite because $W_i$ is
nonzero for only finitely many integers. Therefore the $V_{-j}$ give
a $\zz$-grading on $\Hom_K(S/m^n S, K)$ as an $R$-module.
\end{proof}

We are now ready to present the main result. The
method of the proof will be to show that if $M$ is an $m$-coprimary
module containing an element $u\in 0_{M}^{*}\setminus 0_{M}^{+\GR}$,
then $M$ can be mapped to a finitely generated graded $m$-coprimary
$R$-module where the image of $u$ is not in the plus closure of 0 in
this new module.

\begin{thm}\label{fltc=+cl} Let $R$ be a standard graded
$K$-algebra of characteristic $p>0$. Suppose that $R$ is reduced
(respectively, a domain), and $K$ is perfect (resp., algebraically closed).
If $N_M^* = N_M^F$ (resp., $N_M^* = N_M^+ =
N_M^{+\gr}$) for all finitely generated graded $R$-modules
$N\subseteq
M$ such that $M/N$ is $m$-coprimary, then the same is true for all
finitely generated modules $N\subseteq M$ such that $M/N$ is $m$-coprimary.
\end{thm}
\begin{proof} Let $S=R^\infty$ (resp., $S=R^{+\GR}$). It suffices to
show that $0_M^*\subseteq 0_M^F$ (resp., $0_M^*\subseteq 0_M^{+\gr}$)
when $M$ is $m$-coprimary. Suppose that $u\in 0_M^*\setminus
0_M^F$ (resp., $u\in 0_M^*\setminus 0_M^{+\gr}$). Since $u\not\in
0_M^F$ (resp., $u\not\in 0_M^{+\gr}$), $1\tensor u \neq 0$  via the 
natural map 
$\phi:M\ra S\tensor M$. This implies that there is a surjection of 
the cyclic S-module
$S(1\tensor u)$ onto $K$ sending
$1\tensor u$ to $1\in K$, since the residue field of $S$ is $K$.

Since $\Hom_K(S, K)$ is an injective $S$-module and since we have
a homomorphism $K\ra \Hom_K(S,K)$ that sends 1 to the functional 
that takes $s\in S$ to $s$ modulo $m_S$, there exists a homomorphism
$\psi$ as in the diagram below
$$
\xymatrix{M\ar[r]^{\phi} & S\tensor M\ar@{-->}[rrd]^{\psi} \\
Ru\ar@{^{(}->}[u]\ar[r] & S(1\tensor u)\ar@{^{(}->}[u]\ar[r] & K\ar[r] &
\Hom_K(S,K)
}
$$
such that $\psi\circ\phi(u)\neq 0$. Since $M$ is $m$-coprimary, there exists an $n$
such that $m^n M = 0$. Hence, the image of $M$ under $\psi\circ\phi$
lies in the annihilator of $m^n$ in $\Hom_K(S,
K)$, which is isomorphic to $\Hom_K(S/m^n S, K)$.

By Proposition \ref{grading}, $\Hom_K(S/m^n S, K)$ is a
$\zz$-graded $R$-module. Let $M'$ be the $R$-submodule of
$\Hom_K(S/m^n S, K)$ generated by the homogeneous components of the
images of the generators of $M$. Then $M'$ is a finitely
generated $m$-coprimary graded $R$-module.

Let $\tilde{u}=\psi\circ\phi(u)$, which we
know is nonzero and in $M'$. As $u\in
0_M^*$, we also have that $\tilde{u}\in 0_{M'}^*$. By our hypothesis,
$\tilde{u}\in 0_{M'}^F$ (resp., $\tilde{u}\in 0_{M'}^{+\gr}$)
since $M'$ is graded and $m$-coprimary.
Therefore,
$1\tensor\tilde{u} = 0$ in $S\tensor M'$. Since $\Hom_K(S/m^n S, K)$
is an $S$-module, the inclusion map $M'\incl \Hom_K(S/m^n S, K)$
factors through the map $M' \ra S\tensor M'$, by the universal
property of base change. Thus, the fact that $1\tensor \tilde{u} = 0$ in
$S\tensor M'$ implies that the image of $\tilde{u}$ is 0 in
$\Hom_K(S/m^n S, K)$, a contradiction.
\end{proof}

\section{An Application in Dimension 2}
\label{TCconseq}

In \cite{Br05} and \cite{Br3} (Theorems \ref{Br-ellcurve} and \ref{Br-curve} here),
Brenner shows cases where
the tight closure of a primary homogeneous ideal is the same as
its graded-plus closure. Brenner has observed in correspondence that
it is straightforward to generalize these results
to include finitely generated $m$-coprimary $R$-modules.
The argument is lengthy, like the one for ideals, but the changes are
routine. (The main idea is to replace
the syzygy bundle constructed from homogeneous
generators of an $m$-primary ideal with a syzygy bundle
constructed from homogeneous generators of an $m$-coprimary
submodule $N$ of a graded module $M$. Once one has made
the necessary alterations to \cite[Section 3]{Br03},
all of the relevant proofs in \cite{Br05} and \cite{Br3} follow
seamlessly as they only rely on the aforementioned results
and theorems whose hypotheses only require locally
free sheaves of arbitrary rank, which we obtain in both the ideal
and module cases.)
We state this generalization
as the following theorem.

\begin{thm}[H.\ Brenner] \label{genbrenner} Let $R$ be a
positive characteristic ring. Further, let $R$ be the
homogeneous coordinate
ring of an elliptic curve over an algebraically closed field
$K$, or
let $R$ be any 2-dimensional standard graded $K$-algebra,
where $K$ is the algebraic closure of a finite field. Let
$N\subseteq M$ be finitely generated graded $R$-modules such that
$M/N$ is $m$-coprimary, where $m$ is the homogeneous
maximal ideal of $R$. Then $N_M^* = N_M^+ = N_M^{+\gr}$.
\end{thm}

This result together with Theorem \ref{fltc=+cl} yields an
extension of Theorem \ref{genbrenner}.

\begin{cor} With $R$ as above, $N_M^* = N_M^+ = N_M^{+\gr}$ for all
finitely generated $R$-modules such that $M/N$ is $m$-coprimary.
\end{cor}

Further, if $\Proj\ R$ is an elliptic curve with Hasse invariant 0,
then Brenner's Theorem \ref{Br-ellFr} says that the tight closure of a primary
homogeneous ideal is the same as its Frobenius closure. For example,
this is the case for the cubical cone $R=K[x,y,z]/(x^3+y^3+z^3)$,
when the characteristic of $K$ is congruent to $2\pmod{3}$ (as implied by
\cite[Proposition 4.21]{Ha}). Again,
Brenner's result can be generalized to include finitely generated
homogeneous modules $N\subseteq M$
with $m$-coprimary quotients. This fact can then be
paired with Theorem \ref{fltc=+cl} to give:

\begin{cor} If $R$ is the homogeneous coordinate ring of an
elliptic curve of positive characteristic $p$ with Hasse invariant 0
defined over an algebraically closed field, then $N_M^* = N_M^F$ for
all finitely generated $R$-modules such that $M/N$ is $m$-coprimary.
\end{cor}

For a Noetherian ring $R$ with a maximal ideal $m$,
given the equivalence of tight closure and plus closure
(respectively, graded-plus closure or Frobenius closure) in the $m$-coprimary
case, we can present a characterization of tight closure over
$R_m$ in terms of the $m$-adic closure of certain modules
over $R^+$ (resp., $R^{+\GR}$, $R^{+\gr}$, or
$R^{\infty}$). We start with some general lemmas about tight closure.

\begin{lemma}\label{tclosureinter} Let $(R,m)$ be a reduced local ring
of positive characteristic $p$
that has a test element. Let $M$ be a finitely generated
$R$-module. Then $u\in 0_M^*$ if and only if $u\in \bigcap_k (m^k
M)_M^*$. \end{lemma}
\begin{proof} Let $c$ be a test element in $R$. Then $u\in 0_M^*$ if
and only if $c^{1/q}\tensor u = 0$ in $R^{1/q}\tensor M$ for all $q$
by Lemma \ref{tcchar}. This holds if and only if
$$
c^{1/q}\tensor u \in
\bigcap_k m^k(R^{1/q}\tensor M)
$$
for all $q$ since $(R^{1/q},m^{1/q})$ is also
local, $R^{1/q}\tensor M$ is a finitely generated
$R^{1/q}$-module, and the powers of $mR^{1/q}$ are cofinal with the powers of
$m^{1/q}$. Since $m^k(R^{1/q}\tensor M) = \im(R^{1/q}\tensor
m^kM\ra R^{1/q}\tensor M)$, the above occurs if and only if
$$
c^{1/q}\tensor u\in \im(R^{1/q}\tensor m^kM\ra R^{1/q}\tensor M),
$$
for all $k$ and all $q$. Finally, since $c$ is a test element, the
previous holds if and only if $u\in (m^kM)_M^*$ for all $k$.
\end{proof}

\begin{lemma}\label{madic} Let $R$ be a reduced ring, $I$ an ideal, and
$S=R^\infty$ (respectively, $R$ is also a domain and $S=R^+$ or $R$
is also a graded domain and $S=R^{+\GR}$ or $S=R^{+\gr}$).
Then $u\in (I^kM)_M^F$
(resp., $u\in (I^kM)_M^+$ or $u\in (I^kM)_M^{+\gr}$) for all $k$ if
and only if $1\tensor u\in \bigcap_k I^k(S\tensor M)$. \end{lemma}
\begin{proof} By definition, $u\in (I^kM)_M^F$ (resp., $u\in (I^kM)_M^+$ or $u\in
(I^kM)_M^{+\gr}$) if and only if $1\tensor u\in \im(S\tensor I^kM\ra
S\tensor M)$. This holds if and only if $1\tensor u\in I^k(S\tensor
M)$.
Therefore, $u\in (I^kM)_M^F$ (resp., $u\in (I^kM)_M^+$ or $u\in
(I^kM)_M^{+\gr}$) for all $k$ if and only if $1\tensor u\in \bigcap_k
I^k(S\tensor M)$.
\end{proof}

We now give the promised result connecting tight closure in $R_m$ and
the $m$-adic closure in $R^+$ (resp., $R^{+\GR}$, $R^{+\gr}$, or
$R^{\infty}$).

\begin{prop}\label{tclmadic} Let $R$ be a reduced ring of
characteristic $p>0$. Let $m$ be a maximal ideal of $R$ such that $R_m$
has a test element (e.g., $R_m$ is excellent). Let $S=R^\infty$
(resp., let $R$ also be a domain and $S=R^+$ or let $R$ be a graded
domain and $S=R^{+\GR}$ or $S=R^{+\gr}$).
Moreover, let $R$ be such that Frobenius closure (resp., plus closure
or graded-plus closure) equals tight closure for finitely generated
modules with $m$-coprimary quotient.

Then for any finitely generated
$N\subseteq M$ and $u\in M$, we have $u/1\in (N_m)_{M_m}^*$ if and only if
$1\tensor\overline{u}$ is in the $m$-adic closure of $0$ in $S\tensor M/N$.
For $M$ free, we further note that $(N_m)_{M_m}^*\cap M = \bigcap_k
(N+m^kM)S\cap M$.
\end{prop}
\begin{proof} Since $x\in (N_m)_{M_m}^*$ if and only if $\overline{x}\in
0_{M_m/N_m}^*$ and $M_m/N_m\isom (M/N)_m$,
it is enough to show this for the case $N=0$.
By Lemma \ref{tclosureinter}, $u/1\in 0_{M_m}^*$ if and
only if $u/1\in \bigcap_k (m^k M_m)_{M_m}^*$. Since $M/m^kM$ is clearly
$m$-coprimary, \cite[Proposition 8.9]{HH90} shows that the contraction of
$(m^k M_m)_{M_m}^*$ to $M$ is just $(m^kM)_M^*$ for all $k$. Hence
$u/1\in 0_{M_m}^*$ if and only if $u\in \bigcap_k (m^kM)_M^*$.
By our hypothesis, this holds if and
only if $u\in \bigcap_k (m^k M)_M^F$ (resp., $u\in \bigcap_k (m^k M)_M^+$
or $u\in \bigcap_k (m^k M)_M^{+\gr}$). Then Lemma \ref{madic} shows this
is equivalent to $1\tensor u\in \bigcap_k m^k(S\tensor M)$.

In the case that $M$ is free, the above shows $u\in (N_m)_{M_m}^*\cap M$
if and only if $1\tensor \overline{u}\in \bigcap_k m^k(S\tensor M/N)$, but
$m^k(S\tensor M/N) \isom m^k(MS/NS)$ in this case. Further,
$\overline{u}\in m^k(MS/NS)$ if and only if $u\in (N+m^kM)S$.
\end{proof}

\section{Injective Hulls Over $R^\infty$ and $R^+$}
\label{injhullsect}

In this final section we study the injective hull of the residue field
of $R^\infty$, $R^{+}$, and $R^{+\GR}$,
where $R$ has positive characteristic and
is a complete local domain or
a standard graded $K$-algebra domain.
We start by studying the injective hull, $E_{R^\infty}(K^\infty)$,
where
$R=K[[x_1,\ldots,x_n]]$ or $R=K[x_1,\ldots,x_n]$.
In dimension $n\geq 2$, we show
that there are elements of $E_{R^\infty}(K^\infty)$
that are not killed by any power of the maximal ideal
of $R$.
We then show that this result also holds for complete local domains
and standard graded $K$-algebra domains
in positive characteristic, and that it
holds for the injective hull of the residue field
over $R^+$ or $R^{+\GR}$ as well.
This latter result shows that we cannot
extend Theorem \ref{fltc=+cl} by making use of the injective hull of
the residue field of $R^{+\GR}$ or $R^\infty$ in the analogous way.

\subsection{The Regular Case}

In order to study the injective hull $E_{R^\infty}(K^\infty)$, where
$R=K[[x_{1},\ldots,x_{n}]]$ or $R=K[x_{1},\ldots,x_{n}]$,
we construct a submodule of formal sums
such that the support has DCC. The supports will be subsets of
$-\nn[1/p]^n$, the set of $n$-tuples of nonpositive
rational numbers whose denominators are powers of $p$.
Throughout
the rest of this section, we will use bold letters to stand for $n$-tuples
of elements. We will place a partial
ordering on $n$-tuples by comparing coordinate-wise,
e.g., $\ba > \bb$ if and only if $a_i\geq b_i$, for all $i$,
and $a_j>b_j$, for some $j$. We will define
addition and subtraction of $n$-tuples
as usual. If $\ba\in\qq^n$, then
$\bx^\ba := x_1^{a_1}\cdots x_n^{a_n}$.

\begin{defn} \label{essnotation}
Let $R=K[[\bx]]$ or $R=K[\bx]$, where $K$ is
a field of positive characteristic $p$ and $\dim R=n$.
Given a formal
sum $f = \sum_\ba c_\ba \bx^{-\ba}$, where $\ba\in \nn[1/p]^n$
and $c_\ba\in K^\infty$, we will say that the \textit{support} of $f$
is the subset of $(-\nn[1/p])^n$ given by
$$
\supp(f):= \{ -\ba \, |\, c_\ba\neq 0\}.
$$

Using the same notation, we define the following set of formal sums
$$
N := \{f=\sum_\ba
c_\ba \bx^{-\ba} \, |\, \ba\in\nn[1/p],\,
c_\ba\in K^\infty,\,
\textup{and}\, \supp(f)\, \textup{has DCC} \}.
$$
\end{defn}

\begin{lemma} Using the notation of Definition \ref{essnotation},
$N$ is an $R^\infty$-module with formally defined multiplication.
\end{lemma}
\begin{proof} Let $f_1$, $f_2$ be in $N$, and let $\supp(f_i) = A_i$.
Then
$\supp(f_1+f_2)\subseteq A_1\cup A_2$. Since the union of two
sets with DCC has DCC and a subset of a set with DCC also has DCC,
$f_1+f_2$ is in $N$.
Now, let $s\in R^\infty$. Then
$s\in R^{1/q}$, for some $q=p^e$, so that we can write
$$
s = \sum_{\bb\geq 0}
d_\bb \bx^{\bb/q},
$$
where $\bb\in\nn^{n}$ and $d_\bb\in K^\infty$. Put
$$
f:=\sum_\ba
c_\ba \bx^{-\ba}\in N.
$$
Using formal
multiplication, the coefficient
of $\bx^{-\bs}$ in $sf$ is
$$
\sum_{-\ba+\bb/q = -\bs}
c_\ba d_\bb. \eqno\deqno\label{coeff}
$$
Notice that the coefficient of
 $\bx^{-\bs}$ is 0 if $-\bs = -\ba+\bb/q >0$
as $K^\infty=S/m_S$.

When $-\bs
\leq 0$, for $sf$ to be well-defined, the summation (\ref{coeff})
must consist of
a finite sum of nonzero elements. In the polynomial case, this is
clear. Otherwise, suppose that we have
enumerated the terms contributing to the coefficient
of $\bx^{-\bs}$ and that the set
$$
\{k\in \nn \, |\, -\ba^{(k)} + \bb^{(k)}/q = -\bs,
\text{ and }
c_\ba d_\bb \neq 0\}
$$
is infinite. If there are only finitely many distinct
$\bb^{(k)}$, then (\ref{coeff}) is clearly a finite sum.
We may then assume that there are infinitely
many distinct $\bb^{(k)}$  and thus
assume that all of the $\bb^{(k)}/q$
are distinct. Hence, we obtain an
infinite chain of equalities
$$
-\ba^{(1)} + \bb^{(1)}/q = -\ba^{(2)} + \bb^{(2)}/q =
-\ba^{(3)} + \bb^{(3)}/q = \cdots.
$$
Since the sets $\nn/q$ and $\supp(f)$ have DCC,
we may apply Lemma \ref{infchains} to obtain a
contradiction. Therefore, (\ref{coeff})
is a finite sum, and $sf$ is well-defined.

We also need to show that $\supp(sf)$ has DCC. Suppose to the
contrary that
$$
-\ba^{(1)} + \bb^{(1)}/q > -\ba^{(2)} + \bb^{(2)}/q >
-\ba^{(3)} + \bb^{(3)}/q > \cdots.
$$
is an infinite chain in $\supp(sf)$.
If there are only finitely many distinct
$n$-tuples $\bb^{(k)}/q$, then we also obtain
an infinite descending chain in the
$-\ba^{(k)}$, for $k\gg 0$, a contradiction since
$\supp(f)$ has DCC. We may then assume that
there are infinitely many $\bb^{(k)}/q$ and all are distinct and then apply
Lemma \ref{infchains} again to obtain a contradiction.
\end{proof}

\begin{lemma}\label{infchains} Let $A$ and $B$ be subsets
of $G^n$, where $(G,+)$ is a linearly ordered abelian group. Suppose
that $A$ has DCC and that $B$ has DCC in each coordinate.
If $\{\ba^{(k)}\}_k$ is a sequence
of $n$-tuples in $A$ and $\{\bb^{(k)}\}_k$ is a sequence
of infinitely many distinct $n$-tuples in $B$,
then we cannot obtain an infinite chain
$$
\ba^{(1)} + \bb^{(1)} \geq \ba^{(2)} + \bb^{(2)} \geq
\ba^{(3)} + \bb^{(3)} \geq \cdots. \eqno\deqno\label{ineqchain}
$$
\end{lemma}
\begin{proof} Suppose we have an infinite chain as in (\ref{ineqchain}).
Because each $\bb^{(k)}$ has only finitely
many coordinates and each $\bb^{(k)}$ is distinct, 
after taking subsequences, we may assume
without loss of generality that, for each $i$, either
$b^{(k)}_i = b^{(k')}_i$, or
$b^{(k)}_i < b^{(k+1)}_i$, for all $k,k'$. (The latter assumption
may be made when there are infinitely many distinct
values because $B$ has DCC in each coordinate.)
These conditions imply that
$\bb^{(1)}< \bb^{(2)}< \bb^{(3)} <\cdots$.
If subtract this chain of inequalities
from (\ref{ineqchain}), we obtain
an infinite descending chain
$\ba^{(1)} > \ba^{(2)} > \ba^{(3)} >\cdots$,
which contradicts the fact that $A$ has DCC.
\end{proof}

\begin{prop} \label{Nessext}
Using the notation of Definition \ref{essnotation},
$N$ is an essential extension of $K^\infty$. Therefore,
$N\subseteq E_{R^\infty}(K^\infty)$.
\end{prop}
\begin{proof} The second claim follows immediately from the first.
For the first,
let $f=\sum_{\ba} c_{\ba}
\bx^{-\ba}\in N$.
Since $\supp(f)$ has DCC, we can choose a minimal
element $-\ba^{(0)}$. Then
$\bx^{\ba^{(0)}} \in S$, and
$$
\bx^{\ba^{(0)}} f =
\sum_{\ba} c_{\ba}
\bx^{\ba^{(0)}-\ba}
= c_{\ba^{(0)}} \in K^\infty\setminus \{0\}
$$
as
$\ba^{(0)}_i >\ba_i$, for some $i$, for all
$\ba\neq \ba^{(0)}$ in $\supp(f)$.
\end{proof}

\begin{remk}
M.\ McDermott showed in \cite[Proposition 5.1.1]{McD} that
$N$ is the entire injective hull of
$K^{\infty}$ over $R^{\infty}$ in dimension 1. McDermott's
proof covers the case $R=K[x]$, but the case $R=K[[x]]$ follows
routinely. Whether the result is true for $\dim R\geq 2$ is unknown.
\end{remk}

\begin{prop}\label{badhullreg}
With the notation of Definition \ref{essnotation}, if
 $n\geq 2$, then the injective hull $E_{R^\infty}(K^\infty)$ contains
 an element not killed by any power of
 $m_R = (\bx)R$.
 \end{prop}
 \begin{proof} Let $f = \sum_e x_{1}^{-1/p^e}x_{2}^{-e}$. For $e< e'$,
$-1/p^e < -1/p^{e'}$ and
 $-e > -e'$ so that all elements in $\supp(f)$ are incomparable.
Hence, all chains in $\supp(f)$
 have only one link, and $f\in N$, which injects into
 $E_{R^{\infty}}(K^{\infty})$ by the last proposition.
Now, let $t>0$. Then $x_{2}^t f = \sum_e x_{1}^{-1/p^e}x_{2}^{n-e}$,
 and if $e_0\geq t$, then $t-e_0\leq 0$. Therefore, $x_{2}^t f\neq 0$, and
$m_R^t f\neq 0$, for any $t>0$.
\end{proof}

\subsection{The General Case}

We will now show how we can extend Proposition
\ref{badhullreg} to include complete local domains
and standard graded $K$-algebra domains in
positive characteristic. Moreover, we will also extend
the result to one concerning the injective hull of
the residue field over $R^+$ or $R^{+\GR}$.

An injection of $R$-modules $N\ra M$ is called \textit{pure} if $W\tensor N\ra
W\tensor M$ is an injection for all $R$-modules $W$. When $M/N$ is
finitely presented, the map is pure if and only if the map splits; see
\cite[Corollary 5.2]{HR}. When $S$ is an $R$-algebra and $R\ra S$
is pure as a map of $R$-modules, one calls $S$ \textit{pure} over
$R$.

\begin{lemma}\label{limpure} Let $R=\varinjlim_\alpha R_\alpha$, and let
$S=\varinjlim_\alpha S_\alpha$ such that each $S_\alpha$
is pure over $R_\alpha$. Then $S$ is pure over $R$.
\end{lemma}

If $A$ is a regular ring of positive characteristic
and $R$ is a reduced module-finite extension of $A$, then $A$ is
a direct summand of $R$ as an $A$-module; see \cite[Theorem 1]{Ho73}.
Thus, $A^{1/q}$ is a direct summand of $R^{1/q}$, for all $q=p^e$,
and so the last lemma implies that $R^\infty$ is
pure over $A^\infty$.

\begin{thm}\label{badperfhull} Let $(R,m,K)$ be a complete local
domain (resp., a standard graded $K$-algebra domain) of
positive characteristic and Krull dimension $n\geq 2$.
Then there exists an element of
$E:=E_{R^\infty}(K^\infty)$ that is not killed by any power
of $m$.
\end{thm}
\begin{proof} By the Cohen structure theorem,
$R$ is a module-finite extension of a formal power
series ring $A=K[[x_1,\ldots,x_n]]$ (resp.,
by Noether normalization, $R$ is a module-finite
extension of the graded polynomial ring $A=K[x_1,\ldots,x_n]$).
Since $A$ is regular, $R$ is pure over $A$, and
so the last lemma implies that $R^\infty$ is
pure over $A^\infty$. If we let $E_0:=E_{A^\infty}(K^\infty)$,
then
$$
K^\infty \incl E_0\incl M:=R^\infty\tensor_{A^\infty} E_0.
$$
Since $K^\infty$ is an $R^\infty$-module, we can find
an $R^\infty$-submodule $M'$ of $M$
maximal with respect to not intersecting $K^\infty$. Hence,
$M/M'$ is an essential extension of $K^\infty$ as an $R^\infty$-module.
We can then extend $M/M'$ to a maximal essential extension
$E$ of $K^\infty$ over $R^\infty$. Since the inclusion
$K^\infty\ra E$ factors through $E_0$ and since $E_0$
is an essential extension of $K^\infty$ over $A^\infty$,
$E_0$ injects into $E$ as a map of $A^\infty$-modules.
Since $E_0$ contains an element not killed by any power of
the maximal ideal (resp., the homogeneous maximal
ideal) $m_A$ of $A$ by Proposition \ref{badhullreg}, so does
$E$. Since $m_A$ is primary to $m$, the same element of $E$
not killed by a power of $m_A$ is also not killed by a power of
$m$.
\end{proof}

We can also take advantage of the faithful flatness of
$A^+$ or $A^{+\GR}$ over a regular ring $A$
(see \cite[p. 77]{HH92}) to prove
the existence of elements not killed by a power
of the maximal ideal in the injective hull of the
residue field over $R^+$ or $R^{+\GR}$.

\begin{thm}\label{badinthull} Let $(R,m,K)$ be a complete local
domain (resp., a standard graded $K$-algebra domain) of
positive characteristic and Krull dimension $n\geq 2$.
Then there exists an element of
$E:=E_{R^+}(\overline{K})$ (resp., $E:=E_{R^{+\GR}}(\overline{K})$)
that is not killed by any power of $m$, where $\overline{K}$ is
the algebraic closure of $K$.
\end{thm}
\begin{proof} By the Cohen structure theorem,
$R$ is a module-finite extension of a formal power
series ring $A=K[[x_1,\ldots,x_n]]$ (resp.,
by Noether normalization, $R$ is a module-finite
extension of the graded polynomial ring $A=K[x_1,\ldots,x_n]$).
Thus, $A^+\isom R^+$ (resp., $A^{+\GR}\isom R^{+\GR}$),
and so we may assume that $R=K[[x_1,\ldots,x_n]]$ (resp.,
$R=K[x_1,\ldots,x_n]$). Let $B:=R^+$ (resp., $B:=R^{+\GR}$).
Since $R^{1/q}$ is regular, for all $q$, and since $B$
is a big Cohen-Macaulay $R^{1/q}$-algebra, $B$ is faithfully
flat over $R^{1/q}$. Therefore,
$B$ is flat over $R^\infty$.

The inclusion of $K^\infty\subseteq E_{R^\infty}(K^\infty)$,
together with the flatness of $B$ over $R^\infty$ gives the
following diagram:
$$
\xymatrix{K^\infty\, \ar[d]\ar@{^{(}->}[r] & E_{R^\infty}(K^\infty)\ar[d] \\
B\tensor_{R^\infty} K^\infty\, \ar@{^{(}->}[r] &
B\tensor_{R^\infty} E_{R^\infty}(K^\infty)}
$$
As we have a surjection of $B\tensor_{R^\infty} K^\infty$ onto $\overline{K}$, the
residue field of $B$, we have a map from $B\tensor_{R^\infty} K^\infty$
to $E$, the injective hull of $\overline{K}$ over $B$. Because $E$
is injective, this map lifts to a map from $B\tensor_{R^\infty} E_{R^\infty}(K^\infty)$.
Hence, we obtain a commutative diagram of $R^\infty$-module
maps:
$$
\xymatrix{K^\infty\, \ar@{^{(}->}[r]\ar@{^{(}->}[d] & E_{R^\infty}(K^\infty)\ar[dl] \\
E & }
$$
where the diagonal map is also injective since $E_{R^\infty}(K^\infty)$
is an essential extension of $K^\infty$. Therefore, the element in
$E_{R^\infty}(K^\infty)$ not killed by any power of $m$ (as in
Proposition \ref{badhullreg}) has a nonzero image in
$E$ that is not killed by any power of $m$.
\end{proof}

\section*{Acknowledgments}
I thank Mel Hochster for all of our helpful conversations
during my time at the University of Michigan, and I thank the
anonymous referee for comments that led to clearer expositions
of several proofs.


\begin{thebibliography}{HHH}

\bibitem[AHH]{AHH} \textsc{I. Aberbach}, \textsc{M. Hochster}, and \textsc{C. Huneke},
\textit{Localization of tight closure and modules of finite phantom projective
dimension}, J. Reine. Angew. Math. \textbf{434} (1993), 67--114.

\bibitem[Br1]{Br03} \textsc{H. Brenner}, \textit{Tight closure and projective
bundles}, J. Algebra \textbf{265} (2003), 45--78.

\bibitem[Br2]{Br05} \textsc{H. Brenner}, \textit{Tight closure and plus closure
for cones over elliptic curves}, Nagoya Math. J. \textbf{177} (2005), 31--45.

\bibitem[Br3]{Br3} \textsc{H. Brenner}, \textit{Tight closure and plus closure
in dimension two}, Amer. J. Math. \textbf{128}(2) (2006), 531--539.

\bibitem[Ha]{Ha} \textsc{R. Hartshorne}, \textit{Algebraic Geometry},
Graduate Texts in Mathematics \textbf{52}, Springer-Verlag,
New York, 1977.

\bibitem[Ho]{Ho73} \textsc{M. Hochster}, \textit{Contracted ideals from integral
extensions of regular rings}, Nagoya Math. J. \textbf{51} (1973), 25--43.

\bibitem[HH1]{HH90} \textsc{M. Hochster} and \textsc{C. Huneke}, \textit{Tight
closure, invariant theory, and the Brian\c{c}on-Skoda theorem}, J. Amer.
Math. Soc. \textbf{3} (1990), 31--116.

\bibitem[HH2]{HH92} \textsc{M. Hochster} and \textsc{C. Huneke}, \textit{Infinite
integral extensions and big Cohen-Macaulay algebras}, Annals of Math.
\textbf{135} (1992), 53--89.

\bibitem[HH3]{HH94sm} \textsc{M. Hochster} and \textsc{C. Huneke},
\textit{$F$-regularity, test elements, and smooth base change},
Trans. Amer. Math. Soc. \textbf{346} (1994), no. 1, 1--62.

\bibitem[HH4]{HH95} \textsc{M. Hochster} and \textsc{C. Huneke},
\textit{Applications of the existence of big Cohen-Macaulay algebras},
Adv. Math.  \textbf{113}  (1995),  no. 1, 45--117.

\bibitem[HR]{HR} \textsc{M. Hochster} and \textsc{J. L. Roberts}, \textit{The purity of the
Frobenius and local cohomology}, Adv. Math. \textbf{21} (1976), 117--172.

\bibitem[McD]{McD} \textsc{M. McDermott}, \textit{Tight closure,
plus closure and Frobenius closure in cubical cones}, Thesis,
Univ. of Michigan, 1996.

\bibitem[Sm1]{Sm94} \textsc{K.E. Smith}, \textit{Tight closure of parameter
ideals}, Invent. Math. \textbf{115} (1994), 41--60.

\bibitem[Sm2]{Sm95} \textsc{K.E. Smith}, \textit{Tight closure and graded
integral extensions}, J. Alg. \textbf{175} (1995), 568--574.

\end{thebibliography}
\end{document}